\documentclass[reqno]{amsart}   
\usepackage{pb-diagram,enumerate,amsmath}
 \usepackage{amsfonts,amssymb,amsmath,a4wide}  
\usepackage{pstricks,pst-node,pst-coil}
 
\theoremstyle{plain}    
\newtheorem{thm}{Theorem}[section]   
\newtheorem{cor}[thm]{Corollary}   
\newtheorem{lem}[thm]{Lemma}   
   
\theoremstyle{remark}    
   
\theoremstyle{definition}    
   
\def\al{{\alpha}}       
       
\def\de{{\delta}}

\def\la{{\lambda}}

\def\si{{\sigma}}

\def\ep{{\varepsilon}}

\def\phi{{\varphi}}

\def\fep{f_{\al,\ep}}
\def\gep{g_{\ep}}
\def\insm{\displaystyle\int_M}
\def\vol{dv_g}
\def\volep{dv_{\gep}}
\def\mugep{\mu_1(g_{\ep})}
\def\uep{u_{\ep}}
\def\etar{\eta_{\rho}}
\def\co{c_0}
\def\cep{c_{\ep}}
\def\uepb{\overline{u}_{\ep}}
\def\aep{a_{\ep}}
\def\vep{v_{\ep}}

\DeclareMathAlphabet{\doba}{U}{msb}{m}{n}

\gdef\mR{\doba{R}}       
\gdef\mS{\doba{S}}

\def\lamin{\lambda_{\rm min}^+}   
\def\msup{\mu_{\rm sup}}
\def\Vol{{\mathop{\rm Vol}}}

\long\def\komment#1{}
 
\begin{document}   
\title{The first eigenvalue of Dirac and Laplace operators on surfaces}   
\author{J.F. Grosjean et E. Humbert}   
\date{August 2003} 
\date{\today}   
   
\begin{abstract}   
Let $(M,g,\si)$ be a compact Riemmannian surface equipped with a spin
structure $\si$. For any metric $\tilde{g}$ on $M$, we denote by
$\mu_1(\tilde{g})$ (resp.  $\la_1(\tilde{g})$) the first positive 
eigenvalue of the Laplacian (resp. the Dirac operator) with respect to the metric
$\tilde{g}$. In this paper, we show that 
$$\inf \frac{\la_1(\tilde{g})^2 }{\mu_1(\tilde{g})} \leqslant  \frac{1}{2}.$$
where the infimum is taken over the metrics $\tilde{g}$ conformal to $g$. This answer a question asked by Agricola, Ammann and Friedrich  in \cite{aaf:99}.
\end{abstract}

\maketitle   
\footnote{grosjean@iecn.u-nancy.fr, humbert@iecn.u-nancy.fr} 
   
\tableofcontents   
   
\section{Introduction}   
Let $(M,g,\si)$  be a compact Riemannian  surface equipped with a  spin
structure $\si$. For any metric $\bar{g}$ on $M$, we denote by 
$\Sigma_{\bar{g}} M$ the spinor bundle associated to $\bar{g}$. We let
$\Delta_{\bar{g}} $  be the Laplace-Beltrami operator acting on smooth functions of
$M$ and $D_{\bar{g}}$ be the Dirac operator acting on smooth spinor
fields with respect to  the metric $\bar{g}$. We also
denote by $\mu_1(\bar{g})$ (resp. $\la_1(\bar{g})$)  the smallest
positive eigenvalue of $\Delta_{\bar{g}}$ (resp. $D_{\bar{g}}$). 
Agricola, Ammann  and Friedrich asked
the following question in \cite{aaf:99}: \\

\noindent {\it When $M$ is a  two dimensional torus,
can we find on $M$ a Riemannian metric $\tilde{g}$
for which $\la_1(\tilde{g})^2 < \mu_1(\tilde{g})$ ? }\\

\noindent The main goal of this article is to answer this question. We prove
the 
\begin{thm} \label{main}
There exists a family of  metrics $(g_{\ep})_\ep$ conformal to $g$ for which 
\begin{eqnarray*} \
\limsup_{\ep \to 0} \la_1(g_{\ep})^2  \Vol_{g_{\ep}}(M)
\leqslant     4  \pi
\end{eqnarray*} 
\begin{eqnarray*}
 \liminf_{\ep \to 0} \ \mu_1(g_{\ep})\Vol_{g_{\ep}}(M)  \geqslant   8 \pi.
\end{eqnarray*}
\end{thm}

\noindent Theorem \ref{main} clearly answers the question of
\cite{aaf:99} but says much more: first, the result is true on any
compact Riemannian surface equipped with a spin structure and not only when
$M$ is a two-dimensional torus. In addition, the metric $\tilde{g}$  can be
chosen in a given conformal class. Finally, this metric $\tilde{g}$  can be
chosen such that 
$(2-\de) \la_1(g)^2 < \mu_1(g)$ where $\de>0$ is arbitrary
small. More precisely Theorem \ref{main} shows 

\begin{cor} \label{compar}
On any compact Riemannian surface $(M,g)$, we have 
$$\inf \frac{\la_1(\bar{g})^2 }{\mu_1(\bar{g})} \leqslant  \frac{1}{2}$$
where the infimum is taken over the metric $\bar{g}$ conformal to $g$.
\end{cor}

Theorem \ref{main} has other interesting consequences. Indeed, it proves 
\begin{cor} \label{lamin}
For any compact surface $(M,g)$ equipped with a spin structure $\si$, we
let 
$$\lamin(M,g,\si)= \inf \la_1(\bar{g}) \Vol_{\bar{g}}^\frac{1}{2}(M)$$
where the infimum is taken over the metrics $\bar{g}$ conformal to
$g$. Then, we have $\lamin(M,g,\si) \leqslant  \lamin(\mS^2)$ where
$\lamin(\mS^2)$ is the same invariant computed on the standard sphere $\mS^2$.
\end{cor}
This corollary is an immediate consequence of the fact that $\lamin(\mS^2)= 2
\sqrt{\pi}$ (see \cite{ammann.humbert.morel:p03a}). This result was announced
in \cite{ammann.humbert.morel:p03a}. The conformal
invariant $\lamin$ has been studied in many papers (see for example
\cite{hijazi:86,lott:86,baer:92b,ammann:03,ammann.humbert.morel:p03b, 
ammann.humbert:06} ). Indeed,
  it  
  has many relations with Yamabe problem (see \cite{lee.parker:87}). Corollary
\ref{lamin} has been proved in all dimensions by Ammann in \cite{ammann:03}
if either $n\geqslant 3$ or is $D$ is invertible. Corollary \ref{lamin} extends the result to
the remaining case: $n=2$ and $Ker(D) \not= \{0 \}$. In
\cite{ammann.humbert.morel:p03a}, an alternative proof of the case $n \geq
3$ is given and the proof of the case $n=2$ is skectched.
  
\noindent In the same spirit, a consequence of Theorem \ref{main} is  
\begin{cor} \label{msup}
For any compact surface $(M,g)$, we
let 
$$\msup(M,g)= \sup \mu_1(\bar{g}) \Vol_{\bar{g}}^\frac{1}{2}(M)$$
where the infimum is taken over the metrics $\bar{g}$ conformal to
$g$.  Then, we have $\msup(M,g) \geqslant  \msup(\mS^2)$ where
$\msup(\mS^2)$ is the same invariant computed on the standard sphere $\mS^2$.
\end{cor}
The invariant $\msup$ has been studied in \cite{colbois.elsoufi:03} and
Corollary \ref{msup} is a particular case of Theorem A in this paper. 
We obtain here
another proof.\\

\noindent {\it Acknowledgement:} The authors  are very grateful to Bernd
Ammann for having drawn our attention to the question in \cite{aaf:99}.

\section{Generalized metrics}
Let $f$ be a smooth positive function and set $\bar{g} = f^2 g$. Let also 
for $u \in C^{\infty}(M)$ 
$$I_{\bar{g}}(u) = \frac{\int_M | \nabla u |_{\bar{g}} dv_{\bar{g}} }{\int_M u^2
  dv_{\bar{g}} }.$$ 
It is well known that 
$\mu_1(\bar{g}) = \inf I_{\bar{g}} (u)$ where the infimum is taken over the
smooth non-zero functions $u$ for which $\int_M u dv_{\bar{g}} = 0$. 
We now can write  all these expressions in the metric $g$. We then see that
for  $u \in C^{\infty}(M)$, we have 
$$I_{\bar{g}}(u) = \frac{\int_M |\nabla u|_g^2 dv_g}{\int u^2 f^2 dv_g}$$
and $\mu_1(\bar{g}) = \inf I_{\bar{g}}(u)$ 
where the infimum is taken over the
smooth non-zero functions $u$ for which $\int_M u f^2 dv_g = 0$.
Now if $f$ is only of class $C^{0,a}(M)$ for some $a>0$, we can define  $\bar{g} = f^2 g$. The
2-form $\bar{g}$ is not really a metric since $f$ is not smooth. We then
say that $g$ is a {\it generalized metric}. We can define 
 the first eigenvalue
$\mu_1(\bar{g})$ 
of
$\Delta_{\bar{g}}$ using the definition above. 
Now, by standard methods, one sees that there
exists a function $u \in C^2(M)$ with $\int_M u f^2 dv_g =0$ and such
that $I_{\bar{g}}(u) = \mu_1(\bar{g})$. Writing the Euler equation for $u$,
we see that 
\begin{eqnarray} \label{eq}
\Delta_g u = \mu_1(\bar{g}) f^2 u.
\end{eqnarray}
We prove the following
result
\begin{lem} \label{gene_lapla}
If  $(f_n)$ is  a
sequence of smooth positive functions which converges uniformily to $f$,
then  $\mu_1(f_n^2 g)$ tends to $\mu_1(\bar{g})$.
\end{lem} 
\begin{proof}
Let $u_n$ be a eigenfunction function associated to 
 $\mu_1(f_n^2 g)$. Without loss of
generality, we can assume that $\displaystyle\int_M u_n^2 f_n^2 = 1$. We set $v_n =
 u_n - \frac{\displaystyle\int_M u_n f^2 dv_g}{\displaystyle\int_M f^2 dv_g}$.  We then have 
$\displaystyle\int_M v_n f^2 dv_g = 0$  and hence 
\begin{eqnarray} \label{mu<} 
\mu_1(\bar{g}) \leqslant  I_{\bar{g}} (v_n).
\end{eqnarray}
We have 
$$\int_M |\nabla v_n|^2 dv_ g= \int _M |\nabla u_n|^2 dv_ g = \int_M
u_n \Delta_g u_n dv_g.$$
By equation (\ref{eq}), we get that 
\begin{eqnarray} \label{numer}
\int_M |\nabla v_n|^2 dv_ g= \mu_1(f_n^2 g )  \int_M f_n^2 u_n^2 dv_g=\mu_1(f_n^2 g ).
\end{eqnarray}
We also have 
$$\int_M f^2 v_n^2  dv_g = \int_M f^2 u_n^2 - \frac{{\left( \displaystyle\int_M u_n
      f^2 dv_g \right)}^2 }{\displaystyle\int_M f^2 dv_g}.$$ 
Now,
\begin{eqnarray*}
\left|\int_M u_n
      f^2 dv_g\right|   = \left| \int_M u_n (f^2 - f_n^2)dv_g\right|
\leqslant  C \int_M |u_n| (f + f_n )^2 \parallel f- f_n \parallel_{\infty}. 
\end{eqnarray*}
Since the sequence $(f_n)_n$ tends uniformly to $f$ and since $\displaystyle\int_M f_n^2
u_n^2 dv_g=1 $, we get that 
$\lim_n \displaystyle\int_M u_n
      f^2 dv_g = 0$. In the same way, 
$$\int_M f^2 u_n^2  dv_g = \int_M f_n^2  u_n^2  dv_g+o(1) = 1 + o(1).$$
Finally, we obtain
$$\int_M f^2 v_n^2  dv_g =  1 + o(1).$$
Together with (\ref{mu<}) and (\ref{numer}), we obtain that 
$\mu_1(\bar{g}) \leqslant  \liminf_n \mu_1(f_n^2 g)$.
Now, let $u$ be associated to $\mu_1(\bar{g})$ and set $v= u - \frac{\displaystyle\int_M
  u f_n^2 dv_g}{\displaystyle\int_M f_n dv_g}$. We have $\displaystyle\int_M v^2 f_n^2 dv_g = 0$ and
hence 
$$\mu_1(f_n^2 g) \leqslant  I_{f_n^2 g}(v).$$
It is easy to see that $\lim_n I_{f_n^2 g}(v) = I_{\bar{g}}(u)=
\mu_1(\bar{g})$. We then obtain that  $\mu_1(\bar{g}) \geqslant  \limsup_n
\mu_1(f_n^2 g)$. This proves Lemma \ref{gene_lapla}.
\end{proof}

\noindent In the same way, if $\bar{g} = f^2 g$ is a  metric conformal to $g$
where $f$ is positive and smooth, we define 
$$J'_{\bar{g}}(\psi) = \frac{\displaystyle\int_M |D_{\bar{g}}\psi|_{\bar{g}}^2 f^{-1}
  dv_{\bar{g}}}{\displaystyle\int_M \langle  D_{\bar{g}} 
  \psi, \psi\rangle _{\bar{g}} dv_{\bar{g}}}.$$
The first eigenvalue of the Dirac operator $D_{\bar{g}}$ is then given by 
$\lambda_1^+(\bar{g})  = \inf J'_{\bar{g}}(\psi)$ where the infimum is taken over the
smooth spinor fields $\psi$ for which $  \int_M \langle  D_g \psi, \psi\rangle 
dv_g > 0$. Now, it is well known (see \cite{hitchin:74,hijazi:01})
that where can identify isometrically on each fiber spinor fields for the
metric $g$ and spinor fields for the metric $\bar{g}$. Moreover, we have
for all smooth spinor field:
$$D_{\bar{g}}( f^{-\frac{1}2} \psi )= f^{ -\frac{3}2}D_g \psi.$$
This implies that if we set $\phi = f^{-\frac{1}{2}} \psi$, we have   
$$J_{\bar{g}}(\phi): = \frac{\displaystyle\int_M |D_g\phi|^2 f^{-1} dv_g}{\displaystyle\int_M \langle  D_g \phi, \phi\rangle  dv_g}= J'_{\bar{g}}(\psi)$$
and 
the first eigenvalue of the Dirac operator $D_{\bar{g}}$ is given by 
$\lambda_1^+(\bar{g})  = \inf J_{\bar{g}}(\phi)$ where the infimum is taken over the
smooth spinor fields $\phi$ for which $ \displaystyle \int_M \langle  D \phi, \phi\rangle 
dv_g>0$. With the definition above, we can extend the definition of
$\la_1(\bar{g})$ when $\bar{g}$ is a generalised metric. By standard
methods, there exists a spinor fields  $\phi \in C^1(M)$  such that 
 $\lambda_1^+(\bar{g})= J_{\bar{g}}(\phi)$ and such that 
\begin{eqnarray} \label{eq2} 
D_g \phi = \lambda_1^+(\bar{g}) f \phi.
\end{eqnarray} 
We then have a result  similar to Lemma \ref{gene_lapla}:

\begin{lem} \label{gene_dir}
If  $(f_n)$ is  a
sequence of smooth positive functions which converges uniformly to $f$,
then  $\la_1(f_n^2 g)$ tends to $\la_1(\bar{g})$.
\end{lem}

\noindent The proof is similar to the one of Lemma \ref{gene_lapla} and we omit it here.

\section{The metrics $(g_{\al,\ep})_{\al,\ep}$}
In this paragraph, we construct the metrics $(g_{\al,\ep})_{\al,\ep}$ which
will satisfy:

\begin{eqnarray} \label{main_dir} 
\limsup_{\ep \to 0} \la_1(g_{\al,\ep})^2  \Vol_{g_{\al,\ep}}(M)
\leqslant    4 \pi +C(\al)
\end{eqnarray} 
where $C(\al)$ is a positive constant which goes to $0$ with $\al$ and 
\begin{eqnarray}\label{main_lapla}
 \liminf_{\ep \to 0} \mu_1(g_{\al,\ep})\Vol_{g_{\al,\ep}}(M)  \geqslant   8 \pi.
\end{eqnarray}
Clearly this implies 
Theorem \ref{main}.
By Lemmas \ref{gene_lapla} and \ref{gene_dir}, one can assume that
the metrics $(g_{\al,\ep})_{\al,\ep}$ are generalized 
metrics. We just have to define the volume of $M$ for generalized metric by
$\Vol_{f^2 g} (M)= \int_M f^2 dv_g$. 
At first, without loss of generality, we can assume that $g$ is
flat near  a point $p \in M$. Let $\al>0$ be a small number to be fixed
later such that $g$
is flat on $B_p( \al)$.  We set  
for all $x \in M$ and $\ep >0$,
\[ f_{\al,\ep}(x) = \left\{\begin{array}{ccc}
\frac{\ep^2}{\ep^2 + r^2} & \text{if} & r \leqslant   \al\\
\frac{\ep^2}{\ep^2 + \al^2} & \text{if} & r >  \al
\end{array} \right. \]
where $r = d_g(.,p)$. The function $f_{\al,\ep}$ is of class $C^{0,a}$ for all
$a \in ]0,1[$ and is positive on $M$. We then define for all $\ep >0$, 
$g_{\al,\ep} = f_{\al,\ep}^2 g$. 
The 2-forms $(g_{\al,\ep})_{\al,ep}$ will be the desired
generalized metrics.
For these metrics, we have 
$$\Vol_{g_{\al,\ep}}(M) = \int_M f_{\al,\ep}^2dv_g = \int_{B_p(\al)} f_{\al,\ep}^2 dv_g 
+ \int_{M \setminus B_p(\al)} f_{\al,\ep}^2 dv_g.$$
Since $g$ is flat on $B_p(\al)$, we have 
$$\int_{B_p(\al)} f_{\al,\ep}^2 dv_g = \int_0^{2 \pi} \int_0^{\al}
\frac{\ep^4 r}{(\ep^2 + r^2)^2} dr d\Theta.$$
Setting $y= \frac{x}{\ep}$ we obtain:
$$\int_{B_p(\al)} f_{\al,\ep}^2 dv_g  = 2 \pi \ep^2  \int_0^{\frac{\al}{\ep}}
\frac{r}{(1 + r^2)^2} dr =  2 \pi \ep^2 \left( \int_0^{+ \infty} 
\frac{r}{(1 + r^2)^2} dr +o(1) \right) = \pi \ep^2 + o(\ep^2).$$ 
Since  $f_{\al,\ep}^2  \leqslant  \frac{\ep^4}{\al^4}$ on $M \setminus B_p(\al)$,
 we have $\int_{M \setminus B_p(\al)}
f_{\al,\ep}^2 dv_g = o(\ep^2)$. 
We obtain 
\begin{eqnarray} \label{volume}
 \Vol_{g_{\ep}}(M)=\pi \ep^2 + o(\ep^2).  
\end{eqnarray}
In the whole paper, the notation ``$o(.)$'' must be understood as $\ep$
tends to $0$.
\section{Proof of relation (\ref{main_dir})}
Let $f: \mR^2 \to \mR^2$ be defined by 
$f(x) = \frac{2}{1+ |x|^2}$. Let $\psi_0$ be a non-zero parallel spinor
field on $\mR^2$ such that $|\psi_0|^2 = 1$. As in \cite{ammann.humbert.morel:p03a}, we set on $\mR^2$ 
$$\psi(x) = f^{\frac{n}{2}}(x) (1 - x) \cdot \psi_0.$$ 
As easily computed, we have 
on $\mR^2$
\begin{eqnarray} \label{test_spinor} 
 D\psi  = f \psi \; \hbox{ and } \;  |\psi|= f^{\frac{1}{2}}.
\end{eqnarray}
Now, we fixe a small number $\de >0$ such that $g$ is flat on
$B_p(\de)$. Then, we take $\ep \leqslant  \alpha \leqslant  \de$. We will let $\ep$
goes to $0$. We let also 
$\eta$ be a smooth cut-off function defined on $M$ such that 
$0 \leqslant  \eta \leqslant  1$, $\eta(B_p(\de)) = \{1 \}$, $\eta(M \setminus
  B_p(2 \de))= \{0\}$. Identifying $B_p(\de)$ in $M$ with $B_0(\de)$ in
  $\mR^2$, we can define a smooth spinor field on $M$ by $\psi_{\ep} =
  \eta(x) \psi\left( \frac{x}{\ep}\right)$. 
Using (\ref{test_spinor}), we have 
\begin{eqnarray} \label{dpsi}
D_g(\psi_{\ep}) = \nabla \eta \cdot \psi\left(\frac{x}{\ep}\right) + \frac{\eta}{\ep}
f\left(\frac{x}{\ep}\right) \psi\left(\frac{x}{\ep}\right).
\end{eqnarray}
Since $\langle  \nabla \eta \cdot \psi(\frac{x}{\ep}),  \psi(\frac{x}{\ep})\rangle  \in
i \mR$ and since $|D_g \psi_{\ep}|^2 \in \mR$, we have 
\begin{eqnarray} \label{intdpsi}
\int_M |D_g \psi_{\ep}|^2 f_{\al,\ep}^{-1} dv_g = I_1 + I_2 
\end{eqnarray} 
where 
$$I_1 = \int_M |\nabla \eta|^2 \left| \psi\left(\frac{x}{\ep}\right)\right|^2 dx \; \hbox{ and }
I_2 = \int_M \frac{\eta^2}{\ep^2} f^2 \left(\frac{x}{\ep}\right) \left|
\psi\left(\frac{x}{\ep}\right)\right|^2 f_{\al,\ep}^{-1} dx.$$
At first, let us deal with $I_1$.
By (\ref{test_spinor}), 
$$I_1 \leqslant  C \int_M f\left(\frac{x}{\ep}\right) f_{\al,\ep}^{-1}dx = C \int_{B_p(\al)}
f \left(\frac{x}{\ep}\right) f_{\al,\ep}^{-1}dx +  C \int_{B_p(2\de) \setminus B_p(\al)}
f \left(\frac{x}{\ep}\right) f_{\al,\ep}^{-1}dx$$
where, as in the following, $C$ denotes a constant independant of $\al$ and
$\ep$. On $B_p(\al)$, $f \left(\frac{x}{\ep}\right) f_{\al,\ep}^{-1} = 2$. Hence,
$$ \int_{B_p(\al)}
f \left(\frac{x}{\ep}\right) f_{\al,\ep}^{-1}dx \leqslant  C \al^2.$$
On $B_p(2\de) \setminus B_p(\al)$, since $\ep \leqslant  \al$, 
$$f \left(\frac{x}{\ep}\right) f_{\al,\ep}^{-1} \leqslant  \frac{4\al^2 }{\ep^2 +
  r^2}=\frac{4\al^2 }{\ep^2(1 + \left(\frac{r}{\ep}\right)^2) } $$
Hence,
\begin{eqnarray*}
 \int_{B_p(2\de) \setminus B_p(\al)}
f \left(\frac{x}{\ep}\right) f_{\al,\ep}^{-1}dx & \leqslant  &\frac{4\al^2}{\ep^2}
\int_0^{2 \pi} \int_\al^{\de}
\frac{ r}{(1 +\left(\frac{r}{\ep}\right)^2 )} dr d\Theta\\
& \leqslant  & 8 \pi \al^2 \int_{\frac{\al}{\ep}}^{\frac{\de}{\ep}}
\frac{ r}{(1 +r^2)} dr\\
& \leqslant  & 8 \pi \al^2  \ln{ \left(\frac{\ep^2+ \de^2}{\ep^2 + \al^2}\right)}.
\end{eqnarray*}
We get 
$$\int_{B_p(2\de) \setminus B_p(\al)}
f \left(\frac{x}{\ep}\right) f_{\al,\ep}^{-1}dx\leqslant  
 C\al^2 \ln{\left(\frac{2\de^2}{\al^2}\right) }.$$ 
Finally, we obtain 
\begin{eqnarray} \label{I1}
I_1 \leqslant  C \al^2 + C \ln{\left(\frac{2\de^2}{\al^2}\right) } = a(\al)
\end{eqnarray}
where $a(\al)$ goes to $0$ with $\al$.
Now, by (\ref{test_spinor}),
$$I_2 \leqslant  C \int_{B_p(2\de)}  f^3 (\frac{x}{\ep})f_{\al,\ep}^{-1} dx.$$
Since $f_{\al,\ep} \geqslant   \frac{1}{2} f(\frac{x}{\ep})$, we have 
$$I_2 \leqslant  \frac{2}{\ep^2} \int_{B_p(2\de)} f^2 (\frac{x}{\ep}) dx.$$
Mimicking what we did to get (\ref{volume}), we obtain that 
$$I_2 \leqslant  8 \pi + o(1)$$
when $\ep$ tends to $0$.
Together with (\ref{intdpsi}) and (\ref{I1}), we obtain  
\begin{eqnarray} \label{intdpsi2}
\int_M |D_g \psi_{\ep}|^2 f_{\al,\ep}^{-1} dv_g \leqslant 8 \pi + a(\al) + o(1).
\end{eqnarray}
In the same way, by (\ref{dpsi}), since $\displaystyle\int_M \langle D_g (\psi_\ep), \psi_\ep\rangle 
dv_g \in \mR$ and since $\langle  \nabla \eta \cdot \psi(\frac{x}{\ep}),  \psi(\frac{x}{\ep})\rangle  \in
i \mR$, we have 
$$\int_M \langle D_g (\psi_\ep), \psi_\ep\rangle 
dv_g  =  \int_M \frac{\eta^2}{\ep}   f\left(\frac{x}{\ep}\right) \left|
\psi\left(\frac{x}{\ep}\right)\right|^2 dv_g .$$
By (\ref{test_spinor}), this gives 
$$\int_M \langle D_g (\psi_\ep), \psi_\ep\rangle 
dv_g  = \int_M \frac{\eta^2}{\ep} f^2\left(\frac{x}{\ep}\right)dv_g.$$
With the computations made above, it follows that 
$$\int_M \langle D_g (\psi_\ep), \psi_\ep\rangle 
dv_g = 4 \pi \ep + o(\ep).$$ 
Together with (\ref{intdpsi2}) and (\ref{volume}), we obtain
$$\la_1(g_{\al,\psi})^2  \Vol_{g_{\al,\psi}}(M) \leqslant  
\left(J_{g_{\al,\psi}}(\psi_{\ep})\right)^2 \Vol_{g_{\al,\psi}}(M)  \leqslant 
  \left( \frac{8 \pi + a(\al) + o(1)}{4 \pi \ep+ o(\ep)} \right)^2 (\pi \ep^2 + o(\ep^2))
= \frac{1}{\ep}\left(4 \pi + a(\al)+o(1) \right).$$
Relation (\ref{main_dir}) immediatly follows.

\section{Proof of relation (\ref{main_lapla})}

First we need the following estimate

\begin{lem}\label{poincare} For any $\ep>0$ and $u\in C_c^{\infty}(B_p(\alpha))$, then

$$\insm u^2\fep^2\vol\leqslant\frac{\ep^2}{8}\insm |\nabla u|^2\vol+\frac{1}{\pi\ep^2}\left( \insm u\fep^2\vol\right)^2.$$

\end{lem}

\begin{proof} Let $\gep=\fep^2 g$. Then $(B_p(\al),\gep)$ is embedded in a canonical sphere of volume $\displaystyle\int_{\mR^2}\left(\frac{\ep^2}{\ep^2+r^2}\right)^2dx=\pi\ep^2$. Then from the Poincar\'e-Sobolev inequality, we have

$$\insm u^2\volep\leqslant\frac{1}{\mu_{1,\ep}}\insm|\nabla^{\ep}u|^2_{\gep}\volep+\frac{1}{V_{\ep}}\left(\insm u\volep\right)^2$$

where $\mu_{1,\ep}=\frac{8}{\ep^2}$ is the first nonzero eigenvalue of the Laplacian on the sphere of volume $V_{\ep}=\pi\ep^2$ and $\nabla^{\ep}u$ denotes the gradient of $u$ with respect to the metric $\gep$. Now since $|\nabla^{\ep}u|^2_{\gep}=\fep^{-2}|\nabla u|_g^2$ and $\volep=\fep^2\vol$, we get the desired result.

\end{proof}

\begin{lem}\label{integrat} For any $u,v\in C^{\infty}(M)$, we have

$$\insm (\Delta u)uv^2\vol=\insm |\nabla(uv)|_g^2\vol-\insm u^2|\nabla v|_g^2\vol.$$

\end{lem}

\begin{proof} The proof is an elementary calculation.

\end{proof}

Because of the relation (\ref{volume}), the inequality (\ref{main_lapla}) is equivalent to the following

\begin{align}\label{estimu}\liminf_{\ep\longrightarrow 0}\ep^2 \mugep\geqslant. 8\end{align}

In order to prove this inequality, we assume that for any $\ep$ small enough, there exists $k$, $0<k<1$ so that 

\begin{align}\label{absurde}\mugep <\frac{8}{\ep^2}k.\end{align}

Let $\uep$  be an eigenfunction associated to $\mugep$. Then $\uep\in C^2(M)$ and $\Delta_{\gep}\uep=\mugep\uep$ where $\Delta_{\gep}$ denotes the Laplacian associated to the metric $\gep$. Since the dimension is $2$, $\Delta_{\gep}=\frac{1}{\fep^2}\Delta$ and

\begin{align}\label{equation}\Delta\uep=\mugep\fep^2\uep.\end{align}

We normalize $\uep$ so that $\|\uep\|_{H_1^2}=1$. Up to a subsequence we
can assume that $\insm |\nabla\uep|^2\vol\longrightarrow l$ and
$\insm\uep^2\vol\longrightarrow l'$ with $l+l'=1$. Since $(\uep)$ is bounded
in $H_1^2$, there exists a subsequence so that $\uep\longrightarrow u$
weakly in $H_1^2$. In the following, all the convergences are up to
subsequence. We sometimes omit to recall this fact.

\begin{lem}\label{limc0} There exists a constant $c_0$ such that $u= c_0$.
\end{lem}

\begin{proof}

Let $\varphi\in C^{\infty}(M)$ and

\begin{align*} \eta_{\rho}:=\left\lbrace \begin{matrix} 1 & \text{on} & B_p(\rho)\\
0 & \text{on} & M\setminus B_p(2\rho)\end{matrix}\right.\end{align*}

satisfying $0\leqslant\etar\leqslant 1$ and  $|\nabla\eta_{\rho}|\leqslant\frac{1}{\rho}$. We have

\begin{align}\label{equality}\insm\langle\nabla u,\nabla\varphi\rangle=\insm\langle\nabla u,\nabla(\etar\varphi)\rangle\vol+\insm\langle\nabla u,\nabla((1-\etar)\varphi)\rangle\vol.\end{align}

Now we have

\begin{align*}\insm\langle\nabla u,\nabla(\etar\varphi)\rangle\vol&=\insm\langle\nabla u,\nabla\etar\rangle\varphi\vol+\insm\langle\nabla u,\nabla\varphi\rangle\etar\vol\\
&\leqslant C\left( \int_{B_p(2\rho)}|\nabla u|^2\vol\right)^{1/2}\left(\int_{B_p(2\rho)}|\nabla\etar|^2\vol\right)^{1/2}\\
&+\left( \int_{B_p(2\rho)}|\nabla u|^2\vol\right)^{1/2}\left( \int_{B_p(2\rho)}|\nabla \varphi|^2\vol\right)^{1/2}.
\end{align*}

The limit of the last term is $0$ when $\rho\longrightarrow 0$. Moreover from the definition of $\etar$ and from the fact that $M$ is a $2$-dimensional locally flat domain, the limit of  $\left(\displaystyle\int_{B_p(2\rho)}|\nabla\etar|^2\vol\right)^{1/2}$ is bounded in a neighborhood of $0$. Then we deduce that

\begin{align}\label{intermed}\insm\langle\nabla u,\nabla(\etar\varphi)\rangle\vol\longrightarrow 0\end{align}

 when $\rho\longrightarrow 0$. On the other hand

\begin{align*}\left|\insm\langle\nabla u,\nabla\left((1-\etar)\varphi\right)\rangle\vol\right|&=\lim_{\varepsilon\longrightarrow 0}\left|\insm\langle\nabla \uep,\nabla\left((1-\etar)\varphi\right)\rangle\vol\right|\\
&=\lim_{\varepsilon\longrightarrow 0}\left|\insm(\Delta\uep)(1-\etar)\varphi\vol\right|\\
&=\lim_{\varepsilon\longrightarrow 0}\left|\mugep\insm\fep^2\uep(1-\etar)\varphi\vol\right|.
\end{align*}

Now from the definition of $\fep$ and from (\ref{absurde}) we get

$$\left|\mugep\insm\fep^2\uep(1-\etar)\varphi\vol\right|\leqslant\frac{8}{\ep^2}kC\ep^4\left(\insm\uep^2\vol\right)^{1/2}\left(\insm(1-\etar)\varphi^2\vol\right)^{1/2}$$

where $C$ is a constant depending on the compact support of $(1-\etar)\varphi$. Then making $\varepsilon \longrightarrow 0$, we deduce that

$$\insm\langle\nabla u,\nabla\left((1-\etar)\varphi\right)\rangle\vol=0.$$                                                

Now, reporting this and (\ref{intermed}) in (\ref{equality}) we obtain that $\insm\langle\nabla u,\nabla\varphi\rangle\vol=0$ and $\Delta u=0$ on $M$ in the sense of distributions. This implies that $u\equiv c_0$ on $M$ for a constant $c_0$.

\end{proof}

\begin{lem}\label{cep} Let $(\cep)_{\ep}$ be a bounded sequence of real
  numbers. Then 

$$\insm\fep^2\uep^2\vol\leqslant O(\ep^2\|\uep-\cep\|_{L^2}^2+\ep^4).$$

\end{lem}

                                               \begin{proof} Let $\eta$ be a $C^{\infty}$ function defined on $M$ so that

\begin{align*} \eta:=\left\lbrace \begin{matrix} 1 & \text{on} & B_p(\al/2)\\
0 & \text{on} & M\setminus B_p(\al)\end{matrix}\right.\end{align*}

satisfying $0\leqslant\eta\leqslant 1$ and  $|\nabla\eta|\leqslant 1$.

\noindent {}From the lemma \ref{poincare}, we have

$$\insm(\uep-\cep)^2\fep^2\eta^2\vol\leqslant\frac{\ep^2}{8}\insm|\nabla((\uep-\cep)\eta)|^2\vol+\frac{1}{\pi\ep^2}\left( \insm  (\uep-\cep)\eta\fep^2\vol \right)^2$$

and applying the lemma \ref{integrat}  to the first term of the right hand side, we get

\begin{align*}&\insm(\uep-\cep)^2\fep^2\eta^2\vol\leqslant\\
&\hspace{1cm}\frac{\ep^2}{8}\insm (\Delta(\uep-\cep))(\uep-\cep)\eta^2\vol+\frac{\ep^2}{8}\insm(\uep-\cep)^2|\nabla\eta|^2\vol+ \frac{1}{\pi\ep^2}\left( \insm  (\uep-\cep)\eta\fep^2\vol \right)^2.\end{align*}                    
                                              
{}From (\ref{equation}) we deduce that

\begin{align*}&\insm(\uep-\cep)^2\fep^2\eta^2\vol\leqslant\\
&\hspace{1cm}\frac{\ep^2}{8}\mugep\insm\uep(\uep-\cep)\eta^2\fep^2\vol+\frac{\ep^2}{8}\|\uep-\cep\|_{L^2}^2+ \frac{1}{\pi\ep^2}\left( \insm  (\uep-\cep)\eta\fep^2\vol \right)^2.\end{align*}	

{\em First case : } assume that $\insm\uep(\uep-\cep)\eta^2\fep^2\vol\geqslant 0$.

The relation (\ref{absurde}) implies

$$\insm(\uep-\cep)^2\fep^2\eta^2\vol\leqslant k\insm\uep(\uep-\cep)\eta^2\fep^2\vol+\frac{\ep^2}{8}\|\uep-\cep\|_{L^2}^2+ \frac{1}{\pi\ep^2}\left( \insm  (\uep-\cep)\eta\fep^2\vol \right)^2.$$

A straightforward computation shows that

\begin{align}\label{outil}&(1-k)\insm\uep^2\fep^2\eta^2\vol+\cep^2\insm\fep^2\eta^2\vol\leqslant\notag\\
&\hspace{2cm}(2-k)\cep\insm\uep\fep^2\eta^2\vol+ \frac{\ep^2}{8}\|\uep-\cep\|_{L^2}^2+ \frac{1}{\pi\ep^2}\left( \insm  (\uep-\cep)\eta\fep^2\vol \right)^2.\end{align}

Now note that

\begin{align*}\insm\uep\fep^2\eta^2\vol&=\insm\uep\fep^2(\eta^2-1)\vol+\insm\uep\fep^2\vol\\
&=\insm\uep\fep^2(\eta^2-1)\vol+\frac{1}{\mugep}\insm\Delta\uep\vol\\
&=\insm\uep\fep^2(\eta^2-1)\vol\\
&\leqslant\int_{M\setminus B_p(\alpha/2)}\uep\fep^2(\eta^2-1)\vol\end{align*}

and from the definition of $\fep$ and $\eta$ and from the fact that $\uep$ is bounded in $L^2$, we deduce that

$$\insm\uep\fep^2\eta^2\vol=O(\ep^4).$$

Since $\cep$ is bounded (\ref{outil}) becomes

\begin{align}\label{outil2}(1-k)\insm\uep^2\fep^2\eta^2\vol+\cep^2\insm\fep^2\eta^2\vol&\leqslant O(\ep^4)+\frac{\ep^2}{8}\|\uep-\cep\|_{L^2}^2+ \frac{1}{\pi\ep^2}\left( \insm  (\uep-\cep)\eta\fep^2\vol \right)^2\notag\\
&\hspace{-4cm}=O(\ep^4)+\frac{\ep^2}{8}\|\uep-\cep\|_{L^2}^2+\frac{1}{\pi\ep^2}\left( \insm\fep^2\uep(\eta-1)\vol+\insm\fep^2\uep\vol-\cep\insm\fep^2\eta\right)^2\notag\\
&\hspace{-4cm}=O(\ep^4)+\frac{\ep^2}{8}\|\uep-\cep\|_{L^2}^2+\frac{1}{\pi\ep^2}\left( \insm\fep^2\uep(\eta-1)\vol-\cep\insm\fep^2\eta\right)^2
\end{align}

where in the last equality we have used the fact that $\insm\fep^2\uep\vol=\frac{1}{\mugep}\insm\Delta\uep\vol=0$.

Using the same arguments as above we see that  $\insm\fep^2\uep(\eta-1)\vol=O(\ep^4)$. Reporting this in (\ref{outil2}) we get

$$(1-k)\insm\uep^2\fep^2\eta^2\vol+\cep^2\insm\fep^2\eta^2\vol\leqslant O(\ep^4)+\frac{\ep^2}{8}\|\uep-\cep\|_{L^2}^2+\frac{O(\ep^4)}{\ep^2}\insm\fep^2\eta\vol+\frac{\cep^2}{\pi\ep^2}\left( \insm\fep^2\eta\vol\right)^2.$$

Now

\begin{align*}\insm\fep^2\eta\vol&=\int_{B_p(\al)}\fep^2\vol=\int_0^{2\pi}\int_0^{\al}\frac{\ep^4 r}{(\ep^2+r^2)^2}dr d\Theta\\
&=2\pi\ep^2\int_0^{\al/\ep}\frac{t}{(1+t^2)^2}dt\\
&\leqslant 2\pi\ep^2\int_0^{+\infty}\frac{t}{(1+t^2)^2}dt\\
&=\pi\ep^2.
\end{align*}

This gives

\begin{align*}(1-k)\insm\uep^2\fep^2\eta^2\vol+\cep^2\insm\fep^2\eta^2\vol&\leqslant O(\ep^4)+\frac{\ep^2}{8}\|\uep-\cep\|_{L^2}^2+\frac{\cep^2}{\pi\ep^2}\left( \insm\fep^2\eta\vol\right)^2\\
&=O(\ep^4)+\frac{\ep^2}{8}\|\uep-\cep\|_{L^2}^2+\cep^2\insm\fep^2\eta\vol.
\end{align*}

Finally we have

\begin{align}\label{finally}(1-k)\insm\uep^2\fep^2\eta^2\vol&\leqslant O(\ep^4)+\frac{\ep^2}{8}\|\uep-\cep\|_{L^2}^2+\cep^2\insm\fep^2(\eta-\eta^2)\vol\notag\\
&\leqslant O(\ep^4)+\frac{\ep^2}{8}\|\uep-\cep\|_{L^2}^2+\cep^2\int_{B_p(\al)\setminus B_p(\al/2)}\fep^2\vol\notag\\
&=O(\ep^4)+\frac{\ep^2}{8}\|\uep-\cep\|_{L^2}^2.\end{align}

{\em Second case :} Assume that $\insm\uep(\uep-\cep)\eta^2\fep^2\vol\leqslant 0$. 

In this case, we have

\begin{align}&\insm\uep^2\fep^2\eta^2\vol-2\cep\insm\uep\fep^2\eta^2\vol+\cep^2\insm\fep^2\eta^2\vol\leqslant\notag\\
&\hspace{6cm} O(\ep^4)+\frac{\ep^2}{8}\|\uep-\cep\|_{L^2}^2+ \frac{1}{\pi\ep^2}\left( \insm  (\uep-\cep)\eta\fep^2\vol \right)^2\end{align}

and we conclude as in the previous case.

Then we have proved that

$$\insm\uep^2\fep^2\eta^2\vol=O(\ep^4+\ep^2\|\uep-\cep\|_{L^2}^2).$$

To finish the proof, we write

$$\insm\uep^2\fep^2\vol=\insm\uep^2\fep^2\eta^2\vol+\insm\uep^2\fep^2(1-\eta^2)\vol$$

and the last term is $O(\ep^4)$ which completes the proof.

\end{proof}

\begin{proof} {\em of  Relation (\ref{estimu}).} First we apply the lemma \ref{cep} to $\cep=\co$ and we see that $\co\neq 0$. Indeed, let us compute the $L^2$-norm of the gradient of $\uep$.

\begin{align*} \insm|\nabla\uep|^2\vol&=\insm(\Delta \uep)\uep\vol\leqslant\frac{8k}{\ep^2}\insm\fep^2\uep^2\vol\\
&=\frac{8k}{\ep^2} O(\ep^2\|\uep-\co\|_{L^2}^2+\ep^4)\\
&=o(1).
\end{align*}
 
Then we deduce that up to a subsequence

$$\insm |\nabla \uep|^2\vol\longrightarrow 0.$$
 
 But we have chosen $\uep$ so that $\|\uep\|_{H_1^2}=1$. Then $\|\uep\|_{L^2}\longrightarrow 1$ and $\co\neq 0$.

Now let us consider $\uepb=\frac{1}{\Vol(M)}\insm\uep\vol$ and $\aep=\|\uep-\uepb\|_{H_1^2}$. Then $\uep\longrightarrow \co$ and $\aep\longrightarrow 0$. It follows that the function $\vep=\frac{\uep-\uepb}{\aep}$ satisfies $\|\vep\|_{H_1^2}=1$ and there exists $v\in H_1^2$ so that $\vep\longrightarrow v$ weakly in $H_1^2$ and strongly in $L^2$.

To prove (\ref{estimu}) we will consider two cases.

{\em First case :} Assume that up to a subsequence $\aep=O(\ep)$. 

We have

\begin{align*} \insm(\Delta\uep)^2\vol&=\mugep^2\insm\fep^4\uep^2\vol\\
&\leqslant\mugep^2\insm\fep^2\uep^2\vol\\
&\leqslant\frac{64k}{\ep^4}O(\ep^2\|\uep-\uepb\|_{L^2}^2+\ep^4)\\
&\leqslant \frac{64k}{\ep^4}O(\ep^2\aep^2+\ep^4)\\
&\leqslant M.
\end{align*}

Then $\|\Delta\uep\|_{L^2}$, $\|\nabla\uep\|_{L^2}$ and $\|\uep\|_{L^2}$
are bounded. It well known that the norms 
$$\parallel v \parallel:= \parallel \Delta v\parallel_{L^2}+
\parallel\nabla v \parallel_{L^2}+\parallel v\parallel_{L^2}$$
and  $\parallel v \parallel_{H_2^2}$ are equivalent (it is a direct
consequence of Bochner formula). Hence, this implies that $(\uep)_\ep$ is
bounded in $H_2^2$ which is embedded in $C^0$. Then $\uep\longrightarrow\co$ uniformily up to a subsequence. Since $\co\neq 0$ it follows that for $\ep$ small enough $\uep$ has a constant sign, which is not possible because $\uep$ is an eigenfunction in the metric $\gep$.

{\em Second case :} Assume that $\ep=\aep o(1)$. In this case we have the

\begin{lem} $\vep\longrightarrow c_1$ in $H_1^2$ where $c_1$ is a constant.

\end{lem}

\begin{proof} The proof is similar to this of lemma \ref{limc0}. Indeed we consider $\varphi\in C^{\infty}(M)$ and the function $\etar$ defined in this previous proof. Then

 \begin{align*}\insm\langle\nabla v,\nabla\varphi\rangle=\insm\langle\nabla v,\nabla(\etar\varphi)\rangle\vol+\insm\langle\nabla v,\nabla((1-\etar)\varphi)\rangle\vol.\end{align*}
 
 By the same arguments we have $\insm\langle\nabla v,\nabla(\etar\varphi)\rangle\vol\longrightarrow 0$ when $\rho\longrightarrow 0$. Moreover

 \begin{align*}\left|\insm\langle\nabla v,\nabla\left((1-\etar)\varphi\right)\rangle\vol\right|&=\lim_{\varepsilon\longrightarrow 0}\left|\insm\langle\nabla \vep,\nabla\left((1-\etar)\varphi\right)\rangle\vol\right|\\
&=\lim_{\varepsilon\longrightarrow 0}\left|\insm(\Delta\vep)(1-\etar)\varphi\vol\right|\\
&=\lim_{\varepsilon\longrightarrow 0}\left|\frac{\mugep}{\aep}\insm\fep^2\vep(1-\etar)\varphi\vol\right|.
\end{align*}

Now $\left|\frac{\mugep}{\aep}\insm\fep^2\vep(1-\etar)\varphi\vol\right|\leqslant\frac{8k}{\aep\ep^2}C\ep^4\|\vep\|_{L^2}\left(\insm(1-\etar)\varphi^2\vol\right)^{1/2}$. Since $\ep=\aep o(1)$, we deduce that $\left|\insm\langle\nabla v,\nabla\left((1-\etar)\varphi\right)\rangle\vol\right|=0$ and then $\insm\langle\nabla v,\nabla\varphi\rangle=0$. Therefore $\Delta v=0$ in sense of distributions and $v=c_1$ on $M$.

\end{proof}

Now let $\cep=\uepb+\aep c_1$. Then $\cep\longrightarrow c_0$. We denotes
by $\mu(g)$ the smallest positive eigenvalue of the Laplacian with respect
to the metric $g$. From the definition of $\aep$ and the definition of
$\mu(g)$, we have 
 
\begin{align}\label{denouement}\aep^2\leqslant 2\left( \insm |\nabla\uep|^2\vol+\insm(\uep-\uepb)^2\vol\right)
&\leqslant 2\left( 1+\frac{1}{\mu(g)}\right) \insm|\nabla\uep|^2\vol\notag\\
&=2\left( 1+\frac{1}{\mu(g)}\right)\insm\Delta\uep \uep\vol\notag\\
&=2\left( 1+\frac{1}{\mu(g)}\right)\mugep\insm\fep^2\uep^2\vol.
\end{align}

Applying lemma \ref{cep} we get

\begin{align*}\insm\fep^2\uep^2\vol&=O(\ep^2\|\uep-\cep\|^2_{L^2}+\ep^4)\\
&=O(\ep^2\|\uep-\uepb-\aep c_1\|_{L^2}^2+\ep^4)\\
&=O\left( \aep^2\ep^2\left\|\frac{\uep-\uepb}{\aep}-c_1\right\|_{L^2}^2+\ep^4\right)\\
&=O( \aep^2\ep^2\|\vep-c_1\|_{L^2}^2+\ep^4)\\
&=O(\ep^4)+o(\aep^2\ep^2).
\end{align*}

Now reporting this in (\ref{denouement}) with the estimate (\ref{absurde}) we find

\begin{align*}\aep^2&\leqslant C\frac{8k}{\ep^2}(O(\ep^4)+o(\aep^2\ep^2))\\
&=O(\ep^2)+\aep^2 o(1).
\end{align*}

But $\ep=\aep o(1)$. Then $\aep^2\leqslant C\aep^2 o(1)$ and for $\ep$ small enough $\aep=0$ and $\uep$ is a constant which is impossible.

\end{proof}

\vspace{1cm}             
Authors' address:             
\nopagebreak 
\vspace{5mm}\\ 
\parskip0ex 
%
%
\vtop{ 
\hsize=10cm\noindent 
\obeylines             
Jean-Fran\c cois Grosjean and Emmanuel Humbert,             
Institut \'Elie Cartan BP 239            
Universit\'e de Nancy 1             
54506 Vandoeuvre-l\`es -Nancy Cedex             
France                         
\vspace{0.5cm}             
             
E-Mail:             
{\tt grosjean@iecn.u-nancy.fr}, {\tt humbert@iecn.u-nancy.fr} } 
                  

\begin{thebibliography}{Yam60}           
 

\bibitem[AAF99]{aaf:99}
I. Agricola, B. ammann and T. Friedrich,
\emph{A comparison of the eigenvalues of the Dirac and Laplace operators on
  a two-dimensional torus},
\newblock{Manuscripta Math.}, \textbf{100} (1999), No 2,  231--258.
        

\bibitem[Amm03]{ammann:03}   
B. Ammann, \emph{A spin-conformal lower bound of the first positive {D}irac   
  eigenvalue}, \newblock{{D}iff. {G}eom. {A}ppl.} \textbf{18} (2003), 21--32.   
\bibitem[AH06]{ammann.humbert:06}         
B. Ammann, E. Humbert, 
\newblock{\em The first conformal Dirac eigenvalue on 2-dimensional tori,}
\newblock{ J. Geom. Phys.}, \textbf{56} (2006), No 4, 623--642. 

\bibitem[AHM03]{ammann.humbert.morel:p03a}    
B. Ammann, E. Humbert, B. Morel,    
\newblock{\em A spinorial analogue of Aubin's inequality,}    
\newblock{Preprint}  


\bibitem[AHM03]{ammann.humbert.morel:p03b}         
B. Ammann, E. Humbert, B. Morel,  
\newblock{\em Mass endomorphism and spinorial Yamabe type problems  
on conformally flat manifolds,} to appear in Comm. Anal. Geom.     

   
\bibitem[B\"ar92]{baer:92b}           
C. B\"ar,           
\newblock {\em Lower eigenvalue estimates for Dirac operators,}           
\newblock {Math. Ann.}, 293,  1992.           
    
  

   
\bibitem[CoES03]{colbois.elsoufi:03} 
B. Colbois and A. El Soufi,
\newblock{\em Extremal eigenvalues of the Laplacian in a conformal class of
  metrics: the 'conformal spectrum'}
\newblock{Ann. Global Anal. Geom.}, \textbf{24} (2003), No 4,  337--349.  

\bibitem[Hij86]{hijazi:86}           
O.~Hijazi,           
  \newblock{\em  A conformal lower bound for the smallest eigenvalue of the
    {D}irac operator and {K}illing Spinors },            
\newblock { Comm. Math. Phys.}, \textbf{104} (1986), 151--162. 
           
\bibitem[Hij91]{hijazi:91}             
O.~Hijazi,             
\newblock {\em Premi\`ere valeur propre de l'op\'erateur de {D}irac et nombre de {Y}amabe},             
\newblock { C. R. Acad. Sci. Paris}, \newblock{\em S\'erie I } \textbf{313}, (1991), 865--868.            
 
\bibitem[Hij01]{hijazi:01}
O.~Hijazi,              
\emph{Spectral properties of the {D}irac operator and geometrical structures},
Ocampo, Hernan (ed.) et al., Geometric methods for quantum field theory. 
Proceedings of the summer school, Villa de Leyva, Colombia, July 12-30, 1999. 
Singapore: World Scientific. 116-169,
2001.
  
\bibitem[Hit74]{hitchin:74}
N.~Hitchin, \emph{Harmonic spinors}, Adv. Math. \textbf{14} (1974), 1--55.

\bibitem[LP87]{lee.parker:87}           
J.~M. Lee and T.~H. Parker,           
\newblock {The Yamabe problem,}           
\newblock {\em Bull. Am. Math. Soc., New Ser.}, \textbf{17} (1987), 37--91.           
 
\bibitem[Lot86]{lott:86}           
J. Lott,           
\newblock {\em Eigenvalue bounds for the Dirac operator,}           
\newblock { Pacific J. of Math.}, \textbf{125} (1986), 117--126.             
            
           
           
           
\end{thebibliography}
\end{document}